\documentclass[11pt,a4paper]{amsart}
\usepackage{amssymb}
\usepackage{mathrsfs}

\headsep=1cm \numberwithin{equation}{section}
\newtheorem{thm}{Theorem}[section]
\newtheorem{cor}[thm]{Corollary}
\newtheorem{lem}[thm]{Lemma}
\newtheorem{prop}[thm]{Proposition}
\theoremstyle{definition}

\theoremstyle{remark}

\numberwithin{equation}{section}

\newcommand\Supp{\operatorname{Supp}}
\newcommand\Ass{\operatorname{Ass}}

\newcommand\Rad{\operatorname{Rad}}
\newcommand\Hom{\operatorname{Hom}}
\newcommand\Ext{\operatorname{Ext}}

\begin{document}
\title [Cofiniteness of local cohomology]{Cofiniteness  of local cohomology modules for ideals of dimension one}
\author{Kamal Bahmanpour, Reza Naghipour$^*$  and Monireh Sedghi}
\address{Department of Mathematics, Ardabil branch Islamic Azad University,
P.O. Box 5614633167, Ardabil, Iran.}
\email{bahmanpour@tabrizu.ac.ir}
\address{Department of Mathematics, University of Tabriz, Tabriz, Iran;
and School of Mathematics, Institute for Research in Fundamental
Sciences (IPM), P.O. Box. 19395-5746, Tehran, Iran.}
\email{naghipour@ipm.ir} \email {naghipour@tabrizu.ac.ir}
\address{Department of Mathematics, Azarbaijan University of Tarbiat Moallem, Tabriz, Iran. }%
\email{sedghi@azaruniv.edu}%
\thanks{ 
This research was in part supported by a grant from IPM (No. 900130057 and 900130061).\\
$^*$Corresponding author: e-mail: {\it naghipour@ipm.ir} (Reza Naghipour)}%
\begin{abstract}
Let $R$ denote a commutative Noetherian (not necessarily local)
ring, $M$ an arbitrary $R$-module and $I$ an ideal of $R$ of dimension one. It is shown that
the $R$-module $\Ext^i_R(R/I,M)$ is finitely generated (resp. weakly Laskerian) for all $i\leq {\rm cd}(I,M)+1$
if and only if the local cohomology module $H^i_I(M)$ is $I$-cofinite (resp. $I$-weakly cofinite) for all $i$.
Also, we show that when $I$ is an arbitrary ideal and $M$ is finitely generated module such that the $R$-module $H^i_I(M)$
is weakly Laskerian for all $i\leq t-1$, then $H^i_I(M)$ is $I$-cofinite for all $i\leq t-1$ and for any minimax submodule
$K$  of  $H^{t}_I(M)$,  the $R$-modules $\Hom_R(R/I, H^{t}_I(M)/K)$ and $\Ext^{1}_R(R/I, H^{t}_I(M)/K)$
are finitely generated, where $t$ is a non-negative integer.  This generalizes the main result of Bahmanpour-Naghipour \cite{BN} and 
Brodmann and Lashgari \cite{BL}.

\end{abstract}
\maketitle
\section{Introduction}
Let $R$ denote a commutative Noetherian ring
(with identity) and $I$ an ideal of $R$. For an $R$-module $M$, the
$i^{th}$ local cohomology module of $M$ with respect to $I$ is
defined as$$H^i_I(M) = \underset{n\geq1} {\varinjlim}\,\,
\text{Ext}^i_R(R/I^n, M).$$ We refer the reader to \cite{Gr1} or
\cite{BS} for more details about local cohomology.

In \cite{Gr2} Grothendieck conjectured that for any ideal $I$ of $R$
and any finitely generated $R$-module $M$, the module $\Hom_R(R/I,H^{i}_{I}(M))$ is finitely generated, but soon
Hartshorne was able to present a counterexample to Grothendieck's
conjecture (see \cite{Ha} for details and proof). However, he
defined an $R$-module $M$ to be $I$-cofinite if $\Supp M\subseteq
V(I)$ and $\Ext^{j}_{R}(R/I,M)$ is finitely generated for all $j$
and asked:\\
{\it For which rings $R$ and ideals $I$ are the modules
$H^{i}_{I}(M)$ $I$-cofinite for all $i$ and all finitely generated
modules $M$ }?\\
With respect to this question, Hartshorne in \cite{Ha} and later
Chiriacescu in \cite{Ch} showed that if $R$ is a complete regular
local ring and $I$ is a prime ideal such that $\dim R/I=1$, then
$H^{i}_{I}(M)$ is $I$-cofinite for any finitely generated $R$-module
$M$ (see \cite[Corollary 7.7]{Ha}).

Also, Delfino and Marley \cite[Theorem 1]{DM} and Yoshida \cite
[Theorem 1.1]{Yo} have eliminated the complete hypothesis entirely.
Finally, more recently Bahmanpour and Naghipour \cite[Theorem 2.6]{BN1},
 with completely different methods and relying primarily on induction, removed the local
condition on the ring.

On the other hand,  in \cite{DM1},  Divaani-Aazar and Mafi introduced the class of
weakly Laskerian modules. An $R$-module $M$ is said to
be a {\it weakly Laskerian}, if the set of associated primes of any quotient module $M$ is finite.
Also, they in \cite {DM2} defined an $R$-module $M$ is $I$-{\it weakly cofinite} if
$\Supp M\subseteq V(I)$ and $\Ext^{j}_{R}(R/I,M)$ is weakly Laskerian for all $j$.

One of the aims of the present paper is to prove some new results
concerning cofiniteness  of local cohomology modules $H^{i}_{I}(M)$
for any $R$-module $M$ and any ideal $I$ of dimension one in a commutative Noetherian ring $R$.
More precisely, as a first main result we prove the following:

\begin{thm}
Let $R$ be a commutative Noetherian ring, $M$ an arbitrary $R$-module and $I$ an ideal of $R$, such that
$\dim R/I=1$. Then the following conditions are equivalent:

{\rm(i)} $\Ext^i_R(R/I,M)$ is finitely generated (rep.  weakly Laskerian) for all $i\leq {\rm cd}(I,M)+1.$

{\rm(ii)} $H^i_I(M)$ is $I$-cofinite (resp. $I$-weakly cofinite) for all $i$.

{\rm(iii)} $\Ext^i_R(R/I,M)$ is finitely generated (resp.  weakly Laskerian) for all $i$.

{\rm(iv)} $\Ext^i_R(N,M)$ is finitely generated (resp.  weakly Laskerian) for all $i\leq {\rm cd}(I,M)+1$ and for any finitely generated
$R$-module $N$ with $\Supp N \subseteq V(I)$.

{\rm(v)} $\Ext^i_R(N,M)$ is finitely generated (resp.  weakly Laskerian) for all $i\leq {\rm cd}(I,M)+1$ and for some finitely generated
$R$-module $N$ with $\Supp N = V(I)$.
\end{thm}
Pursuing this point of view further we derive the following consequence of Theorem 1.1.

\begin{cor}
Let $R$ be a commutative Noetherian ring and let $I, J$ be two ideals of $R$ such that
$J\subseteq \Rad(I)$. Let $M$ be a $J$-cofinite (resp. $J$-weakly cofinite) $R$-module.

${\rm (i)}$ If $\dim R/I=1$, then  $H^i_I(M)$ is $I$-cofinite (resp. $I$-weakly cofinite) for all $i$.

${\rm (ii)}$ If $\dim R/J=1$, then  $H^i_J(M)$ is $I$-cofinite (resp. $I$-weakly cofinite) for all $i$.
\end{cor}

We say that an $R$-module $N$ is  {\it minimax module}, if there is a finitely generated
submodule $L$ of $N$, such that $N/L$ is Artinian (see  \cite{Zo1}).  As the second main result, we shall prove the following:

\begin{thm}
Let $R$ be a commutative Noetherian ring, $I$ an ideal of $R$ and $M$ a
finitely generated $R$-module such that for a non-negative integer $t$, the $R$-modules $H^i_I(M)$
are weakly Laskerian for all $i\leq t$. Then the $R$-modules
$H^0_I(M),\dots,H^t_I(M)$ are $I$-cofinite and for any minimax submodule
$K$  of  $H^{t+1}_I(M)$ and for any finitely generated
$R$-module $L$ with $\Supp L\subseteq V(I)$,  the $R$-modules $\Hom_R(L,H^{t+1}_I(M)/K)$ and $\Ext^{1}_R(L,H^{t+1}_I(M)/K)$
are finitely generated.
\end{thm}

Throughout this paper, $R$ will always be a commutative Noetherian
ring with non-zero identity and $I$ will be an ideal of $R$.
For any ideal $\frak a$ of $R$, we denote
$\{\frak p \in {\rm Spec}\,R:\, \frak p\supseteq \frak a \}$ by
$V(\frak a)$. Also,  the {\it radical} of $\frak{a}$, denoted by $\Rad(\frak{a})$, is defined to
be the set $\{x\in R \,: \, x^n \in \frak{a}$ for some $n \in
\mathbb{N}\}$. For any unexplained notation and terminology we refer
the reader to \cite{BS} and \cite{Mat}.\\


\section{The Results}
The main goal of this section are Theorems 2.4 and 2.9.  The following lemmas are needed in the proof of
the main results.  Before bringing them, let us recall that a class $\mathcal{S}$ of $R$-modules is a {\it Serre subcategory} of
the category of $R$-modules, when it is closed under taking submodules, quotients and extensions.
One can easily check that the subcategories of, finitely generated, minimax, weakly Laskerian, and
Matlis reflexive modules are examples of Serre subcategory.

In the case $\mathcal{S}$ is the Serre subcategory of all finitely generated modules, the following lemmas were 
proved by Dibaei and Yassemi. The proofs given in \cite{DY1, DY2} can be easily carried over to an arbitrary 
Serre subcategory $\mathcal{S}$  of the category of $R$-modules.

\begin{lem} {\rm(cf. \cite[Theorem 2.1]{DY1}.)}
    \label{2.1}
 Let $R$ be a Noetherian ring and $I$ an ideal of $R$. Let $s$  be a
 non-negative integer and let $M$ be an $R$-module such that $\Ext_R^s(R/I,M)\in \mathcal{S}$.
If $\Ext_R^j(R/I,H^i_I(M))\in \mathcal{S}$ for all $i<s$ and all $j\geq0$, then
$\Hom_R(R/I,H^s_I(M))\in \mathcal{S}$.
\end{lem}

\begin{lem}{\rm(cf. \cite[Theorem A]{DY2}.)}
    \label{2.2}
Let $R$ be a Noetherian ring and $I$ an ideal of $R$. Let $s$  be a
 non-negative integer and let $M$ be an $R$-module such that $\Ext_R^{s+1}(R/I,M)\in \mathcal{S}$.
If $\Ext_R^j(R/I,H^i_I(M))\in \mathcal{S}$ for all $i<s$ and all $j\geq0$, then
$\Ext^1_R(R/I,H^s_I(M))\in \mathcal{S}$.
\end{lem}

\begin{lem}
 \label{2.3}
 Let $I$ be an ideal of a Noetherian ring $R$ and $M$  a non-zero
$R$-module, such that $\dim M\leq 1$ and $\Supp M\subseteq V(I)$. Then the
following statements are equivalent:

{\rm(i)} $M$ is $I$-cofinite (resp. $I$-weakly  cofinite).

{\rm(ii)} The $R$-modules $\Hom_R(R/I,M)$ and $\Ext^1_R(R/I,M)$ are finitely generated (resp. weakly Laskerian).
\end{lem}
\proof  See \cite[Proposition 2.6]{BNS1} (resp. \cite[Proposition 2.7]{BNS2}).\qed\\

Now we are prepared to state and prove  the first main theorem of this section. Before we state the
Theorem 2.4, let us recall the important concept {\em cohomological
 dimension} of an $R$-module $N$ with respect to an ideal $\frak a$
 of $R$. For an $R$-module $N$, the {\em cohomological
 dimension of } $N$ {\em with respect to} an ideal $\frak a$ of $R$,
 denoted by $\text{cd}(\frak a, N)$, is defined as
 $$\text{cd}(\frak a, N):= \text{sup}\{i\in \mathbb N_0\mid H^i_{\frak a}(N) \neq 0\}.$$

\begin{thm}
 \label{2.4}
Let $R$ be a Noetherian ring, $M$ an  $R$-module and $I$ an ideal of $R$ such that
$\dim R/I=1$. Then the following conditions are equivalent:

{\rm(i)} $\Ext^i_R(R/I,M)$ is finitely generated  (resp. weakly Laskerian) for all $i\leq {\rm cd}(I,M)+1.$

{\rm(ii)} $H^i_I(M)$ is $I$-cofinite (resp. $I$-weakly cofinite) for all $i$.

{\rm(iii)} $\Ext^i_R(R/I,M)$ is finitely generated  (resp. weakly Laskerian) for all $i$.

{\rm(iv)} $\Ext^i_R(N,M)$ is finitely generated  (resp. weakly Laskerian) for all $i\leq {\rm cd}(I,M)+1$ and for any finitely generated
$R$-module $N$ with $\Supp N \subseteq V(I)$.

{\rm(v)} $\Ext^i_R(N,M)$ is finitely generated  (resp. weakly Laskerian) for all $i\leq {\rm cd}(I,M)+1$ and for some finitely generated
$R$-module $N$ with $\Supp N = V(I)$.

{\rm(vi)} $\Ext^i_R(N,M)$ is finitely generated  (resp. weakly Laskerian) for all $i$ and for any finitely generated
$R$-module $N$ with $\Supp N \subseteq V(I)$.

{\rm(vii)} $\Ext^i_R(N,M)$ is finitely generated  (resp. weakly Laskerian) for all $i$ and for some finitely generated
$R$-module $N$ with $\Supp N = V(I)$.

\end{thm}

\proof  We prove theorem for the weakly Laskerian case and by using same proof, the finitely generated case follows. 

In order to prove ${\rm(i)}\Longrightarrow{\rm(ii)}$ we use induction on $i$. When $i=0$, then the exact sequence
$$0\longrightarrow \Gamma_I(M)\longrightarrow M\longrightarrow M/\Gamma_I(M)\longrightarrow0,$$
induces the exact sequence
$$0\longrightarrow \Hom_R(R/I,\Gamma_I(M))\longrightarrow\Hom_R(R/I, M)\longrightarrow \Hom_R(R/I,M/\Gamma_I(M))$$$$\longrightarrow \Ext^1_R(R/I, \Gamma_I(M)) \longrightarrow \Ext^1_R(R/I, M).$$

 As $\Hom_R(R/I,M/\Gamma_I(M))=0$ and $\Ext_R^{i}(R/I,M)$, for $j=0, 1$,  is weakly Laskerian, it follows that
$\Hom_R(R/I, M)$ and $\Ext_R^{1}(R/I,\Gamma_I(M))$ are weakly Laskerian. It now follows from Lemma 2.3 that
 $\Gamma_I(M)$ is $I$-weakly cofinite.

 Assume, inductively, that $i>0$ and that the result has been proved for $i-1$. Then $H^0_I(M),H^1_I(M), \dots, H^{i-1}_I(M)$
are $I$-weakly cofinite, and so by Lemmas 2.1 and 2.2, the $R$-modules $\Hom_R(R/I,H^i_I(M))$ and $\Ext^1_R(R/I,H^i_I(M))$
are weakly Laskerian. Now, it follows again from  Lemma 2.3  that
 $H^i_I(M)$ is $I$-weakly cofinite. The implication ${\rm(ii)}\Longrightarrow{\rm(iii)}$ follows immediately from \cite[Proposition 3.9]{Me2}, and
 the conclusion ${\rm(iii)}\Longrightarrow{\rm(i)}$ is obviously true. Also the implications ${\rm(vi)}\Longrightarrow{\rm(iv)}$ and
${\rm(iv)}\Longrightarrow{\rm(i)}$ are obviously true. For prove ${\rm(iii)}\Longrightarrow{\rm(vi)}$  see \cite[Lemma 2.2]{DM2}, and
the implications ${\rm(vii)}\Longrightarrow{\rm(v)}$ and ${\rm(vi)}\Longrightarrow{\rm(vii)}$ are obviously true. Finally, in order to
complete the proof, we have to show the implication ${\rm(v)}\Longrightarrow{\rm(iv)}$. To this end, let $L$ be a finitely generated
$R$-module with $\Supp L \subseteq V(I)$ and $N$ a  finitely generated $R$-module  such that $\Supp N = V(I)$.

As $\Supp L\subseteq \Supp N$, according to
Gruson's Theorem
 \cite[Theorem 4.1]{Va}, there exists a chain \[0=L_0\subset L_1\subset\cdots\subset L_k=L,\]
 such that the factors $L_j/L_{j-1}$ are homomorphic images of a direct sum of finitely
 many copies of $N$ (note that we may assume that $N$ is faithful). Now consider the exact sequences
\[0\longrightarrow K\longrightarrow
  N^n\longrightarrow L_1\longrightarrow 0  \ \ \ \ \ \ \ \ \ \ \ \ \ \ \ \ \ \ \      \]

\[0\longrightarrow L_1\longrightarrow
   L_2\longrightarrow L_2/L_1\longrightarrow 0  \ \ \ \ \ \ \ \ \ \ \ \ \ \ \ \ \ \ \      \]
                               $$\vdots$$
                                       \[0\longrightarrow L_{k-1}\longrightarrow
  L_k\longrightarrow L_k/L_{k-1}\longrightarrow 0,  \ \ \ \ \ \ \ \ \ \ \ \ \ \ \ \ \ \ \  \]

for some positive integer $n$.  Now,  from the long exact sequence

 \[\cdots\rightarrow {\rm Ext}_R^{i-1}(L_{j-1},N)\rightarrow
   {\rm Ext}_R^{i}(L_j/L_{j-1},N)\rightarrow {\rm
   Ext}_R^{i}(L_j,N) \rightarrow {\rm Ext}_R^i(L_{j-1},N)\rightarrow\cdots,\]

   and an easy induction on $k$, it suffices  to prove the case when
   $k=1$.

 Thus there is an exact sequence \[0\longrightarrow K\longrightarrow
  N^n\longrightarrow L\longrightarrow 0  \ \ \ \ \ \ \ \ \ \ \ \ \ \ \ \ \ \ \      (\ast)\] for some $n\in\Bbb N$
and some finitely generated $R$-module $K$.

    Now, we use induction on $i$. First, $\Hom_R(L,M)$ is a submodule
     of $\Hom_R(N^n,M)$; hence in view of assumption and
    \cite[Lemma 2.2]{DM2}, $\Ext_R^0(L,M)$ is weakly Laskerian.
     So assume that $i>0$ and that $\Ext_R^j(L^\prime, N)$
     is weakly Laskerian for every finitely
  generated $R$-module $L^\prime$ with $\Supp L^\prime\subseteq \Supp N$ and
  all $j\leq i-1$. Now, the exact sequence $(\ast)$ induces the long
  exact sequence \[\cdots\longrightarrow \Ext_R^{i-1}(K,M)\longrightarrow
   \Ext_R^{i}(L,M)\longrightarrow \Ext_R^{i}(N^n,M)\longrightarrow\cdots,\] so that, by the inductive
   hypothesis, $\Ext_R^{i-1}(K,M)$ is weakly Laskerian. On the
   other hand, according to \cite[Lemma 2.2]{DM2}, $\Ext_R^{i}(N^n,M)\cong \stackrel{n}{\oplus}\Ext_R^{i}(L,M)$
   is  weakly Laskerian.  Thus, it follows from \cite[Lemma 2.2]{DM2}
   that $\Ext_R^{i}(L,M)$ is weakly Laskerian,  the inductive step is complete.  \qed\\

\begin{cor}
Let $R$ be a commutative Noetherian ring and let $I, J$ be two ideals of $R$ such that
$J\subseteq \Rad(I)$. Let $M$ be a $J$-weakly cofinite $R$-module.

${\rm (i)}$ If $\dim R/I=1$, then the $R$-module $H^i_I(M)$ is $I$-weakly cofinite for all $i$.

${\rm (ii)}$ If $\dim R/J=1$, then the $R$-module $H^i_J(M)$ is $I$-weakly cofinite for all $i$.
\end{cor}

\proof In order to show that (i), since $J\subseteq \Rad(I)$, it follows that $\Supp R/I\subseteq \Supp R/J$. On the other hand,
since $M$ is $J$-weakly cofinite it follows from \cite[Lemma 2.8]{DM2}
that $M$ is also $I$-weakly cofinite.  Now as $\dim R/I=1$,  the result follows from Theorem 2.4.

Also, to prove (ii), since $\dim R/J=1$ and $M$ is $J$-weakly cofinite it follows from Theorem 2.4 that
$H^i_J(M)$ is $J$-weakly cofinite for all $i$. Now as $\Supp R/I\subseteq \Supp R/J$ it
follows from \cite[Lemma 2.10]{BNS3}  that $H^i_J(M)$ is $I$-cominimax for all $i$. \qed\\

\begin{cor}
Let $R$ be a commutative Noetherian ring and let $I, J$ be two ideals of $R$ such that
$J\subseteq \Rad(I)$. Let $M$ be a  $J$-cofinite $R$-module.

${\rm (i)}$ If $\dim R/I=1$, then the $R$-module $H^i_I(M)$ is $J$-cofinite for all $i$.

${\rm (ii)}$ If $\dim R/J=1$, then the $R$-module $H^i_J(M)$ is $J$-cofinite  for all $i$.
\end{cor}

\proof The result follows from the proof of Theorem 2.4,   Lemma 2.3 and  \cite[Lemma 1]{Ka1}. \qed\\

Before proving  the next main theorem, we need the following lemma and proposition.

\begin{lem}
    \label{2.7}
Let $R$ be a Noetherian ring and $M$ an $R$-module. Then $M$ is weakly Laskerian
if and only if there exists a finitely generated submodule $N$ of $M$ such that
$\Supp M/N$ is finite.
\end{lem}
\proof See \cite[Theorem 3.3]{Ba1}.\qed\\

 The following proposition, which is a generalization the main result of Brodmann and Lashgari,  
 will serve to shorten the proof of the main theorem.

 \begin{prop}
 \label{2.8}
Let $R$ be a Noetherian ring, $I$ an ideal of $R$ and $M$ a
finitely generated $R$-module such that the $R$-modules $H^i_I(M)$
are weakly Laskerian for all $i\leq t$. Then the $R$-modules
$$H^0_I(M),\dots,H^t_I(M)$$ are $I$-cofinite and the $R$-modules $\Hom_R(R/I,H^{t+1}_I(M))$ and $\Ext^{1}_R(R/I,H^{t+1}_I(M))$
are finitely generated. In particular, the  set $\Ass_R H_I^{t+1}(M)/K$ is finite.
\end{prop}
\proof We use induction on $t$. The case $t=0$ easily follows from
Lemmas 2.1 and 2.2. Let $t\geq 1$ and the case $t-1$ is settled. Then by
inductive hypothesis the $R$-modules $$H^0_I(M),\dots,H^{t-1}_I(M)$$
are $I$-cofinite and the $R$-modules $\Hom_R(R/I,H^{t}_I(M))$
and $\Ext^{1}_R(R/I,H^{t}_I(M))$ are finitely generated. Now
since by assumption the $R$-module $H^t_I(M)$ is weakly Laskerian,
it follows from Lemma 2.7 that there is a finitely generated
submodule $N$ of $H^t_I(M)$ such that $\Supp(H^t_I(M)/N)$ is finite
set, and so $\dim H^t_I(M)/N\leq 1$. Now from the
exact sequence $$0\longrightarrow N \longrightarrow H^t_I(M)
\longrightarrow H^t_I(M)/N \longrightarrow 0,$$ it follows that the
$R$-modules $$\Hom_R(R/I,H^t_I(M)/N)\,\,\,\,\,{\rm
and}\,\,\,\,\,\Ext^1_R(R/I,H^t_I(M)/N),$$ are finitely generated.
Therefore it follows from Lemma 2.3 that the $R$-module
$H^t_I(M)/N$ is $I$-cofinite. Hence it follows from the exact
sequence $$0\longrightarrow N \longrightarrow H^t_I(M)
\longrightarrow H^t_I(M)/N \longrightarrow 0,$$ that the $R$-module
$H^t_I(M)$ is $I$-cofinite. Also, it follows from Lemmas 2.1 and 2.2 that the
$R$-modules ${\rm Hom}_R(R/I,H^{t+1}_I(M))$ and ${\rm
Ext}^{1}_R(R/I,H^{t+1}_I(M))$ are finitely generated. This completes
the induction step. \qed\\

Now, we are ready to state and prove the second main theorem of this paper, which is a generalization the main result of  Bahmanpour and Naghipour \cite{BN}.
\begin{thm}
Let $R$ be a commutative Noetherian ring, $I$ an ideal of $R$ and $M$ a
finitely generated $R$-module such that for a non-negative integer $t$, the $R$-modules $H^i_I(M)$
are weakly Laskerian for all $i\leq t$. Then the $R$-modules
$H^0_I(M),\dots,H^t_I(M)$ are $I$-cofinite and for any minimax submodule
$K$ of $H^{t+1}_I(M)$ and for any finitely generated
$R$-module $L$ with $\Supp L\subseteq V(I)$,  the $R$-modules $\Hom_R(L, H^{t+1}_I(M)/K)$ and $\Ext^{1}_R(L, H^{t+1}_I(M)/K)$
are finitely generated.
\end{thm}
\proof By virtue of Proposition 2.8 the $R$-module $H^i_I(M)$ is $I$-cofinite for all $i\leq t$ and the $R$-module
$\Hom_R(R/I,H^{t+1}_I(M))$ is finitely generated. Hence the $R$-module
$\Hom_R(R/I, K)$ is also finitely generated, and so in view of \cite[Proposition 4.3]{Me2}, $K$ is $I$-cofinite.
Thus, \cite[Lemma 1]{Ka1} implies that  $\Ext_R^{i}(L,K)$ is  finitely generated for all $i$.

Next, the exact sequence
\[0\longrightarrow {K}\longrightarrow
   H_I^t(M)\longrightarrow H_I^{t+1}(M)/K \longrightarrow 0\]
provides the following exact sequence,

\[{\Hom_R(L,{\rm H}_I^{t+1}(M))}\longrightarrow
   \Hom_R(L, H_I^{t+1}(M)/K)\longrightarrow  \Ext_R^{1}(L,K)$$ $$\longrightarrow \Ext^1_R(L, H^{t+1}(M))\longrightarrow \Ext^1_R(L, H^{t+1}(M)/K)\longrightarrow    \Ext_R^{2}(L,K) .\]

Now,  since  $\Ext_R^{i}(L, K)$ is finitely generated, so  it is thus
   sufficient for us to show that the $R$-modules  $\Hom_R(L, H_I^{t+1}(M))$ and $\Ext^1_R(L, H^{t+1}(M))$
are finitely generated.  This follows from  Proposition 2.8 and \cite[Lemma 1]{Ka1}. \qed\\
\begin{center}
{\bf Acknowledgments}
\end{center} The authors would like to thank Professor Hossein Zakeri
for his reading of the first draft and valuable discussions. Also, 
we would like to thank from School of Mathematics, Institute for Research in Fundamental
Sciences (IPM), for its financial support.


\end{document}